\documentclass[11pt,a4paper]{scrartcl}
\usepackage{amssymb,amsthm,amsmath,listings,flexisym,breqn}
\theoremstyle{plain}

\newtheorem*{Theorem*}{Theorem}
\newtheorem*{Lemma*}{Lemma}
\theoremstyle{definition}
\newtheorem*{Definition*}{Definition}
\newtheorem*{Remark*}{Remark}
\newtheorem*{Example*}{Example}

\newcommand{\abs}[1]{\lvert#1\rvert}
\DeclareMathOperator{\Gal}{Gal}
\newcommand{\CCC}{\mathbb C}
\newcommand{\FFF}{\mathbb F}            
\newcommand{\QQQ}{\mathbb Q}
\newcommand{\PP}{\mathbb P^1}           
\newcommand{\M}[2]{\mathrm{M}_{#1#2}}             

\begin{document}
\title{A one--parameter family of polynomials with Galois group
  $\pmb{\M24}$ over $\mathbf{\pmb{\QQQ}(t)}$}
\author{Peter M\"uller}
\maketitle
\section{Introduction}
Malle and Matzat proved in \cite[III.7.5]{MM} that the Mathieu group
$\M24$ is a regular Galois extension of $\QQQ(t)$. This group was the
last sporadic group which could be shown to be a Galois group over the
rationals, the only open case remains $\M23$.

The $\M24$ Galois extensions $L/\QQQ(t)$ by Malle and Matzat have the
following properties: Let $\QQQ(t,x)$ be a root field of degree $24$
over $\QQQ(t)$. Then the genus of $\QQQ(t,x)$ is $0$. More precisely,
exactly $4$ places $p_\infty$, $p_0$, $p_1$ and $p_2$ of $\CCC(t)$ are
ramified in $\CCC L$. Let $P_i$ be a place in $\CCC L$ above
$p_i$. Then the inertia group of $P_\infty$ is generated by an element
of order $12$ with two cycles (in the natural degree $24$ action), and
each inertia group of $P_i$, $i=1,2,3$, is generated by an involution
with $8$ fixed points. Clearly, $p_\infty$ is stable under the
absolute Galois group of $\QQQ$, so $p_\infty$ is rational.

We call a $\M24$ Galois extension of $\QQQ(t)$ (or of $\CCC(t)$) with
this ramification data an extension of type $(12,2,2,2)$.

Malle and Matzat show the following: Up to a natural equivalence, the
Galois extensions with the above data, and the added property that
$p_0$ is rational too, are parametrized (up to finitely many
exceptions) by the rational curve $\PP(\QQQ)$. In his thesis
\cite{Granboulan} Granboulan succeeded to compute an explicit
polynomial with these data.

If one requires all branch points $p_i$ to be rational, then the
parametrizing curve has genus $1$, see \cite[III.7.5]{MM}. It was left
open whether this curve has sufficiently many points to give $\M24$
realizations. As a corollary to our computations, we explicitly
compute this curve and show that it is an elliptic curve with positive
rank. So there are infinitely many $\M24$ extension of type
$(12,2,2,2)$ with all branch points rational.

The polynomial by Granboulan is somewhat complicated. We believe that
there are two reasons for this: First, requiring that $p_0$ is
rational adds a further condition, which rules out most polynomials
which could have a nicer shape. Secondly, his polynomial is just a
specialization of an unknown one--parametric family, so the kind of
random specialization could produce ugly coefficients.

For this reason we computed the whole one--parameter family, and
dropped the rationality of $p_0$.

\begin{Theorem*}
Let $t$ be a transcendental over $\QQQ$. For $1\ne s\in\QQQ$ let
$A_s,B_s\in\QQQ[X]$ be as in the Appendix. Then the Galois group of
$(t-A_s(X))^2+(X^2+1)B_s(X)^2$ over $\QQQ(t)$ is $\M24$.
\end{Theorem*}

\begin{Remark*}
For $s=0$ we get the reasonably sized polynomials
\begin{dmath*}
A_0=4194304X^{12} - 72351744X^{10} + 1572864X^9 + 154443776X^8 -
34062336X^7 + 46684160X^6 + 16098816X^5 - 156060348X^4 +
30667728X^3 - 5330757X^2 - 3462498X + 9958791
\end{dmath*}
and
\begin{dmath*}
B_0=-25165824X^{10} - 1572864X^9 + 145227776X^8 - 16515072X^7\\ -
164757504X^6 + 48453120X^5 - 56207872X^4 - 6865152X^3 +
71415384X^2 - 8906760X + 224829.
\end{dmath*}
In the Theorem, we indeed need to exclude the value $s=1$. One can
show that in this case the Galois group is not doubly transitive, so
it is a proper subgroup of $\M24$.
\end{Remark*}
\section{Motivation}
Let $\QQQ(t,x)/\QQQ(t)$ be a $(12,2,2,2)$ extension with branch points
$p_\infty,p_0,p_1,p_2$ as in the introduction, such that the normal
closure has Galois group $\M24$. Since $p_\infty$ is rational, we may
assume that $p_\infty(t)=\infty$. Let $k$ be the residue field of
$P_\infty$. Then $[k:\QQQ]\le 2$ (actually one can show that $k$ is
not real). Then $k(t,x)=k(y)$ (since $\QQQ(t,x)$ has genus $0$ and
$k(t,x)$ has a $k$--rational place). Without loss assume that
$y\mapsto0$ and $y\mapsto\infty$ are the two places above $p_\infty$.

Then
\[
t=\frac{g(y)}{y^{12}}=f(y)
\]
for $g(Y)\in k[Y]$ of degree $24$.  Let $a\mapsto\bar{a}$ be the
automorphism of $kL$ which is the identity on $L$ and has order $2$ on
$k$. From $t=f(y)$ we obtain $t=\bar{t}=\bar{f}(\bar{y})$. On the
other hand, $k(\bar{y})=k(y)$, so
\[
\bar y=\frac{ay+b}{cy+d}
\]
for some $a,b,c,d\in k$. Comparing poles in
$\bar{f}(\frac{ay+b}{cy+d})=f(y)$ shows that either $c=d=0$ or
$a=d=0$. One can show that the former case cannot hold. Thus
$\bar{f}(b/y)=f(y)$. From $y=\bar{\bar
  y}=\bar{b}/\bar{y}=\bar{b}\frac{y}{b}$ we get $\bar{b}=b$, so
$b\in\QQQ$. If $g(Y)=Y^{24}+g_{23}Y^{23}+\dots+g_1Y+g_0$, then
$b^{12}=g_0$. Thus upon replacing $g(Y)$ with $g(\sqrt{b}Y)/b^{12}$,
we may assume that $g(Y)$ is monic with $g_0=1$, for the price that we
possibly need to replace $k$ by a quadratic extension which we still
call $k$. Note also that $\overline{g_{12}}=g_{12}$.

The ramification over $p_0,p_1,p_2$ translates to
\begin{equation}\label{E:AB}
g(Y)-\omega_iY^{12}=A_i(Y)^2B_i(Y),
\end{equation}
where $D(T)=(T-\omega_0)(T-\omega_1)(T-\omega_2)\in k[T]$, and
$A_i,B_i\in k[\omega_i][Y]$ are of degree $8$.

Generically, the branch points $\omega=\omega_0, \omega_1, \omega_2$
are conjugate over $k$. Let $\sigma\in\Gal(\bar\QQQ/k)$ with
$\omega^{\sigma^i}=\omega_i$. Write $A_0=U_0+\omega U_1+\omega^2U_2$
with $U_0,U_1,U_2\in k[Y]$, and similarly $B_0=V_0+\omega
V_1+\omega^2V_2$. Then $A_i=U_0+\omega_i U_1+\omega_i^2U_2$ and
$B_i=V_0+\omega_i V_1+\omega_i^2V_2$. By adding a constant to
$g(Y)/Y^{12}$ we may assume that
$D(T)=(T-\omega_0)(T-\omega_1)(T-\omega_2)=T^3+pT+q$.

We set up a system of $49$ equations for $50$ unknowns: The
polynomials $U_1,U_2,V_1$ and $V_2$ have degree $7$, while $U_0$ and
$V_0$ are monic of degree $8$. This gives $48$ unknown
coefficients. The remaining $2$ unknowns are $p$ and $q$ in
$D(T)=T^3+pT+q$.

The $49$ equations come from noting that \eqref{E:AB} is equivalent to
$A_i(Y)^2B_i(Y)+\omega_iY^{12}$ being independent of $i$, that is
$(U_0+\omega U_1+\omega^2U_2)^2(V_0+\omega V_1+\omega^2V_2)$ having
all coefficients in $k$ (and constant term $1$). Thus treating
$\omega$ as a variable and reducing the expansion of $(U_0+\omega
U_1+\omega^2U_2)^2(V_0+\omega V_1+\omega^2V_2)+\omega_Y^{12}$ modulo
$\omega^3+p\omega+q$ gives $49$ equations.

One can revert the arguments which led to the special shape of $f(Y)$,
similar arguments appear in \cite{Granboulan}. Suppose that $\bar
f(Y)=f(\frac{b}{Y})$ with $b\in\QQQ$. Let $y$ be a root of
$f_s(Y)-t$. Set $x=(y+\bar y)/2$, $z=(y-\bar y)/2$, so $y=x+z$ with
$\bar x=x$ and $\bar z=-z$. From $f(y)=t$ and $\bar t=t$ we get
$f(\frac{b}{\bar y})=t$. But $y$ is the single root of $f(Y)-t$ in
$k(y)$, so $\frac{b}{\bar y}=y$. We get $x^2+z^2=y\bar y=b$. Next
write
\[
t=f(y)=\frac{g(x+z)}{(x+z)^{12}}=g(x+z)\frac{(x-z)^{12}}{b^{12}}.
\]
Expand the right hand side, and reduce modulo $x^2+z^2-b$ with respect
to $z$. This yields $t=A(x)+zB(x)$ for $A,B\in\QQQ[X]$. Then
\[
0=(t-A(x))^2-z^2B(x)^2=(t-A(x))^2+(x^2-b)B(x)^2,
\]
so $(t-A(X))^2+(X^2-b)B(X)^2$ is a minimal polynomial for $x$ over
$\QQQ(t)$.
\section{The Computation}
Granboulan gives a single rational function $g(Y)/Y^{12}$ which after
adding a constant has the shape as above, with $b=-1$. The fact that
all functions $g(Y)/Y^{12}$ with these data and monodromy group $\M24$
are parametrized by a rational curve means that the $49$ equations in
$\CCC^{50}$ (we cannot fix the at most biquadratic field $k$, since it
could vary with $g$) describe a rational curve. Furthermore, the
coefficients $1=\gamma_0,\gamma_1,\dots,\gamma_{23},\gamma_{24}=1$ of
$g(Y)$ are polynomials in the unknowns we work with. So if we fix a
pair $1\le i<j\le 23$, then the coefficients $\gamma_i$ and $\gamma_j$
should be related by a genus $0$ curve equation. Using a
straightforward Newton iteration and starting from Granboulan's
example, we moved in small steps and computed examples of $g(X)$ for
$\gamma_{12}=25,26,27,\dots,124$ to a very high precision ($3000$
binary bits). Fix an index $1\le i\le 22$. Let
$(\gamma_{23,j},\gamma_{i,j})$ be the $100$ pairs of the corresponding
coefficients. We tried if there is a polynomial relation of total
degree $m$ by minimizing $\abs{\sum_{r+s\le
    m}a_{r,s}\gamma_{23,j}^r\gamma_{i,j}^s}$ subject to $a_{0,0}=1$
for unknowns $a_{r,s}\in\CCC$. For the first $m$ which gave a good
approximation we used the function \texttt{algdep} from the computer
algebra package Sage \cite{sage} to determine if $a_{r,s}$ can be
expected to be algebraic. To our surprise, it turned out that all the
coefficients actually lie in $\QQQ(i)$ (so $k$ from above is quadratic
and does not vary with $g$), and that each $\gamma_i$ is a polynomial
in the imaginary part of $\gamma_{23}$. The family $g_s(Y)$,
$s\in\CCC$, which we obtained fulfills $\overline{f_s}(Y)=f_s(-1/Y)$
for all real $s$. Now retrieve $A_s$ and $B_s$ from
$f_s(Y)=g_s(Y)/Y^{12}$ as described in the previous section. The
polynomials $A_s$ and $B_s$ from the appendix slightly differ from
those just obtained: The Galois group of $(t-A(X))^2+(1+X^2)B(X)^2$
doesn't change if we multiply $A$ and $B$ by the same nonzero factor
from $\QQQ$, and it does not change if we add an element from $\QQQ$
to $A$.

The monic cubic $D(T)$ whose roots are the finite branch points of the
splitting field of $(t-A_s(X))^2+(1+X^2)B_s(X)^2$ has the form
$D(T)=T^3+p(s)T^2+q(s)T+r(s)$, where $p,q,r\in\QQQ[S]$ have degrees
$24$, $44$, and $68$, respectively. (Since Sage cannot handle
polynomials of large degree, we computed $p$, $q$, and $r$ with the
help of Magma \cite{magma}.) Note that by adding a constant to $A(X)$,
which amounts to adding the same constant to $t$, we have given up the
original condition that the coefficient of $T^2$ in $D(T)$
vanishes. If we want to make a branch point $p_0$ rational, then we
need to find a condition on $s$ such that $D(T)$ has a rational
root. By Malle--Matzat, the $\M24$ extensions with rational $p_0$ are
still parametrized by $\PP(\QQQ)$. Thus the curve $T^3+p(S)T+q(S)=0$
should be a rational curve, so there should be a cubic rational
function $S(Z)\in\QQQ(Z)$ such that $T^3+p(S(Z))T+q(S(Z))=0$ has a
root in $\QQQ(Z)$. Let $T_0$ be a root of $T^3+p(S)T+q(S)$. Working
out the ramification of the cubic extension $\QQQ(S,T_0)/\QQQ(S)$
allows us to express $S$ in terms of $Z$, where
$\QQQ(S,T_0)=\QQQ(Z)$. Indeed, $S(Z)=\frac{{\left(1-2 \, Z\right)}
  {\left(9 \, Z^{2} + 16 \, Z + 21\right)}}{25 \, {\left(Z^{2} +
    1\right)}}$.

Eventually, we want to get the condition that all finite branch points
are rational. The above parametrization gives a linear factor of
$D(T)=T^3+p(S(Z))T+q(S(Z))$. The square--free part of the discriminant
of the quadratic co--factor is $\frac{1}{81}(Z + 2)(81 \, Z^{3}
+ 36 \, Z^{2}\linebreak + 122 \, Z - 2)$. Thus we need to study
rational points on the hyperelliptic genus $1$ curve
$W^2=\frac{1}{81}{\left(Z + 2\right)}\left(81 \, Z^{3} + 36 \, Z^{2} +
122 \, Z - 2\right)$. The substitution $Z=\frac{30}{U}-2$,
$W=\frac{50V}{3U^2}$ puts this curve in the Weierstrass form
$V^2=U^3-38U^2+540U-2916$. The point $(u,v)=(30,78)$ is on this curve,
but it is not a torsion point by the Nagell--Lutz Theorem ($78$ does
not divide the discriminant of $U^3-38U^2+540U-2916$). Thus there are
infinitely many $\M24$ extensions of $\QQQ(t)$ of type $(12,2,2,2)$
and all branch points rational. An example is given by $s=21/25$.

It remains to verify that the polynomials in the Theorem have the
correct Galois group. Since $\M24$ is self--normalizing in $S_{24}$, it
suffices to show that the given polynomials have Galois group $\M24$
over $\CCC(t)$. Identifying the Galois group with the monodromy group,
it is clear that the group does not change if we vary $s$ along a path
such that each $F_s(X)=(t-A(X))^2+(1+X^2)B(X)^2$ has $4$ distinct
branch points. Computing the discriminant of $D(T)$ from above gives a
high degree polynomial in $s$ whose single rational root is $s=1$.

Thus it suffices to show that $F_0(X)$ has the correct monodromy
group. One checks that $F_0(X)$ is irreducible. If we set $t=1$ and
factor over $\FFF_7$, we get a linear factor and an irreducible factor
of degree $23$. By the Dedekind criterion, the Galois group of
$F_0(X)$ contains a $23$--cycle, so it is doubly
transitive. Furthermore, one verifies that the discriminant of $F_0$
is a square in $\QQQ[t]$. So the Galois group of $F_0$ is either
$\M24$ or $A_{24}$. To rule out the latter case, one can numerically
compute the four generators of the monodromy group corresponding to
four branch points. In a forthcoming paper, we develop an algebraic
criterion which bounds the Galois group from above and which is
applicable here.
\section{Appendix}
Here we give the polynomials $A_s$ and $B_s$ from the theorem. These
and further polynomials and discriminants, corresponding to $p_0$
being rational or all $p_i$ rational, appear on the website
\texttt{www.mathematik.uni-wuerzburg.de/\~{}mueller/Papers/m24.sage}\ .

\begin{align*}
A_s &=
4194304 X^{12} + 25165824 s X^{11} + (-66060288 s^{4} - 132120576 s^{3} + 125829120 s^{2}\\
& + 157286400 s - 72351744) X^{10}+ (7864320 s^{6} - 291766272 s^{5} - 530055168 s^{4}\\
&  + 374865920 s^{3} + 589824000 s^{2} - 257163264 s + 1572864) X^{9} + (95068160 s^{8}\\
& + 404160512 s^{7} - 307265536 s^{6} - 1711472640 s^{5} - 319438848 s^{4}\\
& + 1774780416 s^{3} + 300646400 s^{2} - 680329216 s + 154443776) X^{8}\\
& + (-19660800 s^{10} + 245915648 s^{9} + 1056260096 s^{8} + 425984 s^{7} - 2924347392 s^{6}\\
& - 1609285632 s^{5} + 2559344640 s^{4} + 1031438336 s^{3} - 1117782016 s^{2} + 308436992 s\\
& - 34062336) X^{7} + (-24121344 s^{12} - 187533312 s^{11} + 41097216 s^{10} + 1319323648s^{9}\\
& + 1257576448 s^{8} - 1814122496 s^{7} - 3242037248 s^{6} - 111439872 s^{5} + 2190483456 s^{4}\\
& + 603674624 s^{3} - 701696000 s^{2} + 22325248 s + 46684160) X^{6} + (6127104 s^{14}\\
& - 18555136 s^{13} - 233881088 s^{12} - 175877632 s^{11} + 782083584 s^{10} + 1554816256 s^{9}\\
& + 305029632 s^{8} - 2126697472 s^{7} - 2687739392 s^{6} + 342526208 s^{5} + 2244343296 s^{4}\\
& + 129376768 s^{3} - 926989824 s^{2} + 248273664 s + 16098816) X^{5} + (750020 s^{16}\\
& + 14359744 s^{15} + 12053088 s^{14} - 114504896 s^{13} - 207573520 s^{12} - 47353152 s^{11}\\
& + 587791264 s^{10} + 1525759808 s^{9} + 817520152 s^{8} - 2767906752 s^{7} - 2971570528 s^{6}\\
& + 2058861504 s^{5} + 2304495344 s^{4} - 1478264768 s^{3} - 673535648 s^{2} + 730274240 s\\
& - 156060348) X^{4} + (-240000 s^{18} - 771520 s^{17} + 4552208 s^{16} + 5407552 s^{15}\\
& - 12085824 s^{14} - 57342016 s^{13} - 224008704 s^{12} - 145445568 s^{11} + 869959872 s^{10}\\
& + 1405781568 s^{9} - 746917472 s^{8} - 2364491840 s^{7} + 148282176 s^{6} + 1768511296 s^{5}\\
& - 573188480 s^{4} - 739240000 s^{3} + 693032896 s^{2} - 247013248 s + 30667728) X^{3}\\
& + (8267 s^{20} - 60292 s^{19} - 189418 s^{18} + 387764 s^{17} + 49279 s^{16} + 10416944 s^{15}\\
& + 7799240 s^{14} - 96630576 s^{13} - 173543962 s^{12} + 115412424 s^{11} + 495715012 s^{10}\\
& + 346585368 s^{9} - 264331194 s^{8} - 710757904 s^{7} - 246179640 s^{6} + 577521040 s^{5}\\
& + 203985999 s^{4} - 432517988 s^{3} + 97491606 s^{2} + 11254228 s - 5330757) X^{2}\\
& + (336 s^{22} + 2884 s^{21} + 2918 s^{20} - 280 s^{19} - 187260 s^{18} - 551084 s^{17}\\
& + 3439942 s^{16} + 3990720 s^{15} - 13251728 s^{14} - 31779288 s^{13} - 22616724 s^{12}\\
& + 2200016 s^{11} + 154829368 s^{10} + 220854120 s^{9} - 178041252 s^{8} - 448632704 s^{7}\\
& + 189685920 s^{6} + 382850836 s^{5} - 184274706 s^{4} - 241599992 s^{3} + 191501284 s^{2}\\
& - 37543740 s - 3462498) X
\end{align*}
\begin{align*}
B_s &=
(25165824 s^{2} + 25165824 s - 25165824) X^{10} + (-1572864 s^{4} + 132120576 s^{3}\\
& + 116391936 s^{2} - 119537664 s - 1572864) X^{9} + (-99614720 s^{6} - 304349184 s^{5}\\
& + 268959744 s^{4} + 804782080 s^{3} - 132120576 s^{2} - 452198400 s + 145227776) X^{8}\\
& + (16711680 s^{8} - 350552064 s^{7} - 1014300672 s^{6} + 369819648 s^{5} + 2081488896 s^{4}\\
& - 21430272 s^{3} - 1230962688 s^{2} + 417398784 s - 16515072) X^{7} + (59408384 s^{10}\\
& + 340467712 s^{9} - 186449920 s^{8} - 1965195264 s^{7} - 1035206656 s^{6} + 2817441792 s^{5}\\
& + 2004549632 s^{4} - 1769635840 s^{3} - 729972736 s^{2} + 797122560 s - 164757504) X^{6}\\
& + (-13969920 s^{12} + 98482176 s^{11} + 646695936 s^{10} + 245766144 s^{9} - 2084156928 s^{8}\\
& - 2286686208 s^{7} + 2005444608 s^{6} + 2897768448 s^{5} - 1126577664 s^{4} - 1002874880 s^{3}\\
& + 866995200 s^{2} - 278562816 s + 48453120) X^{5} + (-5992448 s^{14} - 66636032 s^{13}\\
& - 15972864 s^{12} + 521017856 s^{11} + 692872192 s^{10} - 646882048 s^{9} - 1987540480 s^{8}\\
& - 873563136 s^{7} + 1382451200 s^{6} + 1778852096 s^{5} - 371893760 s^{4} - 916518400 s^{3}\\
& + 335021056 s^{2} + 142912256 s - 56207872) X^{4} + (1621440 s^{16} + 424640 s^{15}\\
& - 47240768 s^{14} - 50228288 s^{13} + 157415360 s^{12} + 409560000 s^{11} + 274859712 s^{10}\\
& - 417022784 s^{9} - 1449295296 s^{8} - 722824128 s^{7} + 1625973056 s^{6} + 1180120896 s^{5}\\
& - 1145317056 s^{4} - 494068416 s^{3} + 737746496 s^{2} - 163911616 s - 6865152) X^{3}\\
& + (5000 s^{18} + 1616880 s^{17} + 2639064 s^{16} - 12017344 s^{15} - 21008416 s^{14}\\
& - 17153152 s^{13} + 56032160 s^{12} + 274166336 s^{11} + 340729008 s^{10} - 518906656 s^{9}\\
& - 1093524848 s^{8} + 245942208 s^{7} + 1311410912 s^{6} - 338458624 s^{5} - 772161888 s^{4}\\
& + 347966656 s^{3} + 363724040 s^{2} - 309525584 s + 71415384) X^{2} + (-16656 s^{20}\\
& - 101144 s^{19} + 120440 s^{18} + 90968 s^{17} - 360488 s^{16} - 1152768 s^{15} - 16231616 s^{14}\\
& - 20363328 s^{13} + 68279648 s^{12} + 152960784 s^{11} - 28705040 s^{10} - 273001168 s^{9}\\
& - 97695248 s^{8} + 225560384 s^{7} + 24231552 s^{6} - 162720512 s^{5} + 97096176 s^{4}\\
& + 102722504 s^{3} - 152893224 s^{2} + 65919672 s - 8906760) X + 527 s^{22} + 1138 s^{21}\\
& - 549 s^{20} - 1556 s^{19} - 32451 s^{18} + 372274 s^{17} + 599009 s^{16} - 2859792 s^{15}\\
& - 6690746 s^{14} + 2971140 s^{13} + 16373566 s^{12} + 15848168 s^{11} + 5309674 s^{10}\\
& - 29793004 s^{9} - 62702414 s^{8} + 18214832 s^{7} + 82023963 s^{6} - 11931446 s^{5}\\
& - 56636393 s^{4} + 22914700 s^{3} + 5802345 s^{2} - 4202118 s + 224829
\end{align*}


\end{document}